\documentclass[10pt]{article}

\usepackage{amsmath}
\usepackage{amssymb}

\usepackage{graphicx}

\usepackage{cite}

\usepackage{color} 

\usepackage{verbatim} 

\usepackage[labelfont=bf,labelsep=period,justification=raggedright]{caption}

\date{}

\newcommand{\R}{\mathbb{R}}
\newcommand{\N}{\mathbb{N}}

\newcommand{\Prob}{\mathbb{P}}
\newcommand{\E}{\mathbb{E}}
\newcommand{\Fskript}{\mathcal{F}}

\newenvironment{Proof}[1][]
{\protect\noindent\textbf{Proof#1.} }
{\hfill $\blacksquare$ \\}

\begin{document}

\title{Feller Processes: The Next Generation in Modeling. Brownian Motion, L\'evy Processes and Beyond}

\author{Bj\"orn B\"ottcher\\ \textit{\small TU Dresden, 
Institut f{\"u}r mathematische Stochastik,}\\
\textit{\small 01062 Dresden,
Germany,
{bjoern.boettcher at tu-dresden.de}}}

\date{}

\maketitle

\section*{Abstract}
We present a simple construction method for Feller processes and a framework for the generation of sample paths of Feller processes. The construction is based on state space dependent mixing of L\'evy processes.

Brownian Motion is one of the most frequently used continuous time Markov processes in applications. In recent years also L\'evy processes, of which Brownian Motion is a special case, have become increasingly popular.

L\'evy processes are spatially homogeneous, but empirical data often suggest the use of spatially inhomogeneous processes. Thus it seems necessary to go to the next level of generalization: Feller processes. These include L\'evy processes and in particular Brownian motion as special cases but allow spatial inhomogeneities.

Many properties of Feller processes are known, but proving the very existence is, in general, very technical. Moreover, an applicable framework for the generation of sample paths of a Feller process was missing. We explain, with practitioners in mind, how to overcome both of these obstacles. In particular our simulation technique allows to apply Monte Carlo methods to Feller processes.


\section*{Introduction}

The paper is written especially for practitioners and applied scientists. It is based on two recent papers in stochastic analysis \cite{BoetSchi2009,SchiSchn09}. We will start with a survey of applications of Feller processes. Thereafter we recall some existence and approximation results. In the last part of the introduction we give the necessary definitions. 

The main part of the paper contains a simple existence result for Feller processes and a description of the general simulation scheme. These results will be followed by several examples. 

The source code for the simulations can be found as supporting information (Appendix S1).

\subsection*{Motivation}

Brownian motion and more general L\'evy processes are used as models in many areas: For example in medicine to model the spreading of diseases \cite{Jans1999}, in genetics in connection with the maximal segmental score \cite{DoneMaller2005a}, in biology for the movement patterns of various animals (cf. \cite{ReynFrye2007} and the references therein), for various phenomena in physics \cite{Woyc2001} and in financial mathematics \cite{ContTank2004}. 
In these models the spatial homogeneity is often assumed for simplicity, but empirical data or theoretical considerations suggest that the underlying process is actually state space dependent. Thus Feller processes would serve as more realistic models. We give some explicit examples: 
\begin{itemize}
\item In hydrology stable processes are used as models for the movement of particles in contaminated ground water. It has been shown that state space dependent models provide a better fit to empirical data \cite{Zhan2007,Reev2008}. Also on an intuitive level it seems natural that different kinds of soils have different properties. Thus the movement of a particle should depend on its current position, i.e.\ the soil it is currently in.
\item In geology also stable processes are used in models for the temperature change. Based on ice-core data the temperatures in the last-glacial and Holocene periods are recorded. Statistical analysis showed that the temperature change in the last-glacial periods is stable with index 1.75 and in the Holocene periods it is Gaussian, i.e. stable with index 2 (see Fig. 4 in \cite{Ditl1996}).
\item For a technical example from physics note that the fluctuations of the ion saturation current measured by Langmuir probes in the edge plasma of the Uragan-3M stellarator-torsatron are alpha-stable and the alpha depends on the distance from the plasma boundary \cite{Gonc2003}.
\item Anomalous diffusive behavior has been observed in various physical systems and a standard model for this behavior are continuous time random walks (CTRWs) \cite{Meer2002a,Gore2004}. To study these systems the limiting particle distribution is a major tool, which is in fact a Feller process \cite{Baeu2005}. 
\item In mathematical finance the idea of extending L\'evy processes to L\'evy-like Feller processes was first introduced in \cite{BarnLeve2001}. The proposed procedure is simple: A given L\'evy model usually uses a parameter dependent class of L\'evy processes. Now one makes the parameters of the L\'evy process (in its characteristic exponent) \textit{price}-dependent, i.e.\ the increment of the process shall depend on the current \textit{price}. This procedure is applicable to every class of L\'evy processes, but the existence has to be shown for each class separately \cite{Bass1988a,BarnLeve2001,BoetJaco2003}. 
\end{itemize}

Thus there is plenty of evidence that Feller processes can be used as suitable models for real-world phenomena.

\subsection*{Existence and Approximation}

Up to now general Feller processes were not very popular in applications. This might be due to the fact that the existence and construction of Feller processes is a major problem. There are many approaches: Using the Hille-Yosida theorem and Kolmogorov's construction \cite{Jaco20012005,Hoh98hab}, solving the associated evolution equation (Kolmogorov's backwards equation) \cite{Boet2008,Boet2005,Koch1989,Kolo2000}, proving the well-posedness of the martingale problem \cite{Bass1988a,Hoh98hab,Stro1975}, solving a stochastic differential equation \cite{IkedWata1989,JacoShir2002,Stro2003}.
The conditions for these constructions are usually quite technical. Nevertheless, let us stress that the proof of the very existence is crucial for the use of Feller processes. Some explicit examples to illustrate this will be given at the end of the next section.

Our construction will not yield processes as general as the previous ones, but it will still provide a rich class of examples. In fact the presented method is just a simple consequence of a recent result on the solutions to certain stochastic differential equations \cite{SchiSchn09}.  

Furthermore each of the above mentioned methods also provides an approximation to the constructed Feller process. Most of them are not usable for simulations or work only under technical conditions. Also further general approximation schemes exist, for example the Markov chain approximation in \cite{MaRoeckZhan1998}. But also the latter is not useful for simulations, since the explicit distribution of the increments of the chain is unknown. 

In contrast to these we derived in \cite{BoetSchi2009} a very general approximation scheme for Feller processes which is also usable for simulations. We will present here this method for practitioners.

\subsection*{L\'evy processes and Feller processes}

Within different fields the terms \textit{L\'evy process} and \textit{Feller process} are sometimes used for different objects. Thus we will clarify our notion by giving precise definitions and mentioning some of the common uses of these terms.

A stochastic process is a family of random variables indexed by a time parameter $t\in[0,\infty)$ on a probability space $(\Omega, \Fskript, \Prob)$.  For simplicity we concentrate on one-dimensional processes. The expectation with respect to the measure $\Prob$ will be denoted by $\E$.\\
Although this will not appear explicitly in the sequel, a process will always be equipped with its so-called natural filtration, which is a formal way of taking into account all the information related to the history of the process. Technically the filtration, which is an increasing family of sigma fields indexed by time, is important since a change from the natural filtration to another filtration might alter the properties of the process dramatically.\\

A {\bf L\'evy process} starting in $L_0:=0$ is a stochastic process $(L_t)_{t\in[0,\infty)}$ with
\begin{enumerate}
\item[-] independent increments: The random variables $L_{t_1}, L_{t_2} - L_{t_{1}}, L_{t_3} - L_{t_{2}}, \ldots$ are independent for every increasing sequence $(t_n)_{n\in\N}$, 
\item[-] stationary increments: $L_t-L_s$ has the same distribution as $L_{t-s}$ for all $s < t$, 
\item[-] c\`adl\`ag paths: Almost every sample path is a right continuous function with left limits.
\end{enumerate}
For equivalent definitions and a comprehensive mathematical treatment of L\'evy processes and their properties see \cite{Sato99}.

Note that the term \textit{L\'evy flight} often refers to a process which is a continuous time random walk (CTRW) with spatial increments from a  one-sided or two-sided stable distribution (the former is also called L\'evy distribution). In our notion the processes associated with these increments are L\'evy processes which are called \textit{stable subordinator} and \textit{stable process}, respectively.

A L\'evy process $(L_t)_{t\in[0,\infty)}$ on its probability space is completely characterized by its L\'evy exponent $\xi \mapsto \psi(\xi)$ calculated via the characteristic function
$$\E(e^{i\xi L_t}) = e^{-t\psi(\xi)}.$$
The most popular L\'evy process is Brownian motion ($\psi(\xi) = \frac{1}{2} |\xi|^2$), which has the special property that almost every sample path is continuous. In general, L\'evy processes have discontinuous sample paths, some examples with their corresponding exponents are the Poisson process ($\psi(\xi)=e^{i\xi}-1$), the symmetric $\alpha$-stable process ($\psi(\xi)=|\xi|^\alpha$ with $\alpha \in (0,2]$), the Gamma process ($\psi(\xi)=\ln(1+i\xi)$) and the normal inverse Gaussian process ($\psi(\xi) = -i\mu\xi + \delta\left(\sqrt{\alpha^2-(\beta+i\xi)^2}-\sqrt{\alpha^2-\beta^2}\right)$ with $0\leq|\beta|<\alpha,\ \delta>0,\ \mu\in\R$). 

Classes of L\'evy exponents depend, especially in modeling, on some parameters. Thus one can easily construct a family of L\'evy processes by replacing these parameters by state space dependent functions. Another approach to construct families of L\'evy processes is to introduce a state space dependent mixing of some given L\'evy processes. We will elaborate this in the next section.

Given a family of L\'evy processes $\left((L_t^{(x)})_{t\in[0,\infty)}\right)_{x\in\R}$, i.e.\ given a family of characteristic exponents $(\psi_x)_{x\in\R}$, we can construct for fixed $x_0\in\R$, $T\in[0,\infty)$, $n\in\N$ a Markov chain as follows:
\begin{enumerate}
\item[1.] The chain starts at time 0 in $x_0$.
\item[2.] The first step is at time $\frac{1}{n}$ and it is distributed as $L_\frac{1}{n}^{(x_0)}$. The chain reaches some point $x_1$.
\item[3.] The second step is at time $\frac{2}{n}$ and it is distributed as $L_\frac{1}{n}^{(x_1)}$. The chain reaches some point $x_2$.
\item[4.] The third step is at time $\frac{3}{n}$ and it is distributed as $L_\frac{1}{n}^{(x_2)}$. The chain reaches some point $x_3$.
\item[5.] ... etc. until time $T$.
\end{enumerate}

This Markov chain is spatially inhomogeneous since the distribution of the next step always depends on the current position. 
If the chain converges (in distribution for $n\to \infty$ and every fixed $T\in[0,\infty)$) then the limit is - under very mild conditions (see \cite{Cher1974} and also Theorem 2.5 by \cite{JacoPotr2004}) - a Feller process. Formally, a {\bf Feller process} is a stochastic process $(X_t)_{t\in[0,\infty)}$ such that the operators 
$$T_t f(x):= \E(f(X_t)|X_0 = x), \ \ t\in[0,\infty), x\in \R$$
satisfy 
$$T_0 = id,\ T_s \circ T_t = T_{s+t}\ (s,t \geq 0)\ \ \text{ (semigroup property)}$$
and 
$$ \lim_{t\to 0}\ \sup_{x\in\R}\big|T_t f(x) - f(x)\big| = 0 \ \text{ (strong continuity)}$$
for all $f$ which are continuous and vanish at infinity.
 
A Feller processes is sometimes also called: L\'evy-type process, jump-diffusion, process generated by a pseudo-differential operator, process with a L\'evy generator or process with a L\'evy-type operator as generator. Note that in mathematical finance often the Cox-Ingersoll-Ross process \cite{CoxIngeRoss1985} is called \textit{the Feller process}, but in our notion this is a Feller diffusion in the sense of \cite{Fell1951}. For a comprehensive mathematical treatment of Feller processes and their properties see \cite{Jaco20012005}.

The generator $A$ of a Feller process is defined via
\begin{equation*} 
\lim_{t\to 0}\ \sup_{x\in\R}\left|Af(x)- \frac{T_t f(x)- f(x)}{t}\right| = 0
\end{equation*}
for all $f$ such that the limit exists. Moreover, if the limit exists for arbitrarily often differentiable functions with compact support then the operator $A$ has on these functions the representation
$$ Af(x) = - \int_\R e^{ix\xi}\ \psi_x(\xi)\ \int_\R (2\pi)^{- 1} e^{-iy\xi} f(y)\,dy\ d\xi,$$
where for each fixed $x\in \R$ the function $\xi \mapsto \psi_x(\xi)$ is a L\'evy exponent. Thus a family of L\'evy processes with L\'evy exponents $(\psi_x)_{x\in\R}$ corresponds to the Feller process $(X_t)_{t\in[0,\infty)}$ with generator $A$ as above.

If the corresponding family of L\'evy processes is a subset of a \textit{named} class of L\'evy processes, one calls the Feller process also by the name of the class and adds \textit{-like} or \textit{-type} to it. Thus for example a Feller process corresponding to a class of symmetric stable processes is called symmetric stable-like process.

In general, as mentioned in the previous section, the construction of a Feller process corresponding to a given family of L\'evy processes is very complicated. It even might be impossible as the following examples show: 
Let $(L_t^{(x)})_{t\geq 0}$ be the family of L\'evy processes with characteristic exponents $\psi_x(\xi) = -ia(x)\xi$, i.e. the L\'evy processes have deterministic paths $L_t^{(x)} = a(x)t$. Now if $a(x)=x$ a corresponding Feller process exists, starting in $x$ it has the path $X_t = xe^{t}$. But for $a(x)=x^2$ and $a(x)=\mathrm{sgn}(x)\sqrt{x}$ a corresponding Feller process does not exist, $a(x)=x^2$ yields paths which do not tend to negative infinity as $x\to-\infty$ and $a(x)=\mathrm{sgn}(x)\sqrt{x}$ yields paths which are not continuous with respect to the starting position.

However, we will present in the next section a very simple method to construct Feller processes.  

\section*{Results and Discussion}

\subsection*{Construction of Feller processes by mixing L\'evy processes}

Suppose we know (for example based on an empirical study) that the process we want to model behaves like a L\'evy process $(L^{(1)}_t)_{t\geq 0}$ in a region $K_1$ and like a different L\'evy process $(L^{(2)}_t)_{t\geq 0}$ in a region $K_2$. Then we know that a Feller process which models this behavior exists by the following result:

\noindent{\bf Theorem.}
If the sets $K_1$, $K_2$ are uniformly separated, i.e.\ there exists an $\varepsilon>0$ such that
$$\inf_{x\in K_1,\, y\in K_2} \|x-y\| > \varepsilon$$
then there exists a Feller process $(X_t)_{t\geq 0}$ which behaves like $L^{(1)}$ on $K_1$ and like $L^{(2)}$ on $K_2$.\\
\begin{Proof} Let $\psi^{(i)}$ be the characteristic exponent of $(L_t^{(i)})_{t\geq 0}$ for $i=1,2$.
Under the above condition there exist non-negative bounded and Lipschitz continuous functions $a^{(1)}$ and $a^{(2)}$ such that
$$a^{(1)}(x) = 1 \text{ for all } x\in K_1 \text{ and } a^{(1)}(x) = 0 \text{ for all } x\in K_2,$$  
$$a^{(2)}(x) = 0 \text{ for all } x\in K_1 \text{ and } a^{(2)}(x) = 1 \text{ for all } x\in K_2.$$  

Now set for $\xi_1,\xi_2,x \in \R$ 
$$\psi\left(
\begin{array}{c}
\xi_1\\
 \xi_2	
\end{array}
\right) := \psi^{(1)}(\xi_1)+\psi^{(2)}(\xi_2),\ 
\Phi(x) := \left(\begin{array}{c}
a^{(1)}(x) \\
 a^{(2)}(x)
\end{array}
\right)^\top 
\in \R^{1\times 2}
$$
and note that for $x,\xi \in\R$ 
$$\psi(\Phi^\top(x)\xi) = \psi^{(1)}\left(a^{(1)}(x)\xi\right)+\psi^{(2)}\left(a^{(2)}(x)\xi\right)$$
holds. Thus corresponding to the family of L\'evy processes defined by the L\'evy exponents
$$\psi_x(\xi):=\psi^{(1)}\left(a^{(1)}(x)\xi\right)+\psi^{(2)}\left(a^{(2)}(x)\xi\right)=\psi(\Phi^t(x)\xi)$$
there exists a Feller process as a consequence of Corollary 5.2 from \cite{SchiSchn09} and $\psi_x(\xi) = \psi^{(i)}(\xi)$ for $x\in K_i$ holds ($i=1,2$), i.e.\ $(X_t)_{t\geq 0}$ behaves like $L^{(i)}$ on $K_i$ for $i=1,2$.
\end{Proof}

Note that the theorem extends to any finite number of L\'evy processes $(L_t^{(i)})_{t\geq 0}$ $(i=1,..,n)$ with corresponding regions $K_i$. More generally for any finite number of independent L\'evy processes $(L_t^{(i)})_{t\geq 0}$ $(i=1,..,n)$ with corresponding characteristic exponents $\psi^{(i)}$ and non-negative bounded and Lipschitz continuous functions $x\mapsto\alpha^{(i)}(x)$ the family $(\psi_x)_{x\in\R}$ with
$$\psi_x(\xi) := \psi^{(1)}\left(\alpha^{(1)}(x)\xi\right)+\psi^{(2)}\left(\alpha^{(2)}(x)\xi\right)+\ldots +\psi^{(n)}\left(\alpha^{(n)}(x)\xi\right)$$
defines a family of L\'evy processes $\left((\tilde L_t^{(x)})_{t\geq 0}\right)_{x\in\R}$ and there exists a corresponding Feller process $(X_t)_{t\geq 0}$.

To avoid pathological cases one should assume $\alpha^{(1)}(x) + \alpha^{(2)}(x) + \ldots + \alpha^{(n)}(x) >  0$ for all $x$. Further note that the following equality in distribution holds for all $x$ and $t$
$$\tilde L_t^{(x)} \stackrel{d}{=} \alpha^{(1)}(x) L_t^{(1)} + \alpha^{(2)}(x) L_t^{(2)}+\ldots+\alpha^{(n)}(x) L_t^{(n)}.$$
Thus if one knows how to simulate increments of the $L_t^{(i)}$ one can also simulate increments of $\tilde L_t^{(x)}$. We will see in the next section that simulation of increments of the corresponding family of L\'evy processes is the key to the simulation of the Feller process.

\subsection*{Simulation of Feller processes}
\label{sim_section}

Given a Feller process $(X_t)_{t\in[0,\infty)}$ with corresponding family of L\'evy processes $\left((L_t^{(x)})_{t\in[0,\infty)}\right)_{x\in\R}$ we can use the following scheme to approximate the sample path of $X_t$:

\begin{enumerate}
\item[1.] Select a starting point $x_0$, the time interval $[0,T]$ and the time-step size $h$.
\item[2.] The first point of the sample path is $x= x_0$ at time $t=0$.  
\item[3.] Draw a random number $z$ from the distribution of $L_{h}^{(x)}$ ($x$ is the \emph{current} position of the sample path).
\item[4.] The next point of the sample path is $x \rightsquigarrow x+z$ at time $t \rightsquigarrow t+h$.
\item[5.] Repeat 3. and 4. until $t\in \left(T-h,T\right]$.
\end{enumerate}
  
The simulated path is an approximation of the sample path of the Feller process, in the sense that for $h\to 0$ it converges toward the sample path of the Feller process on $[0,T]$. To be precise, for the convergence the Feller process has to be unique for its generator restricted to the test functions and the family of L\'evy processes $L_{t}^{(x)}$ has to satisfy some mild condition on the $x$-dependence: The L\'evy exponent $\psi_x(\xi)$ has to be bounded by some constant times $1+|\xi|^2$ uniformly in $x$, see \cite{BoetSchi2009} for further details. This condition is satisfied for many common examples of Feller processes, in particular for the processes constructed in the previous section.

The reader familiar with the Euler scheme for Brownian or L\'evy-driven stochastic differential equations (SDEs) will note that the approximation looks like an Euler scheme for an SDE. In fact it is an Euler scheme, but the corresponding SDE does not have such a nice form as for example the L\'evy-driven SDEs discussed in \cite{ProtTala97}. This is due to the fact that in their case for a particular increment all jumps of the driving term are transformed in the same manner, but in the general Feller case the transformation of each jump can depend explicitly on the jump size. More details on the relation of this scheme to an Euler scheme can be found in a forthcoming paper \cite{BoetSchn2009}.

\subsection*{Examples}
\label{examples}

We will now present some examples of Feller processes together with simulations of their sample paths. The first example will show the generality of the mixture approach, the remaining examples are special cases for which the existence has been shown by different techniques. 

All simulations are done with the software package R \cite{R2010} and the source code of the figures can be found as supporting information (Appendix S1).

\subsubsection*{Brownian-Poisson-Cauchy-mixture Feller process}
To show the range of possibilities which are covered by the mixture approach we construct a process which behaves like
\begin{align*}
&\text{ Brownian motion } &\text{on } (-\infty,-6),\\
&\text{ a Poisson process} &\text{on } (-4,4),\\
&\text{ a Cauchy process} &\text{on }  (6,\infty).
\end{align*}
For this we just define a family of L\'evy processes by the family of characteristic exponents $(\psi_x)_{x\in\R}$ with
$$\psi_x(\xi) := a_1(x) \frac{1}{2}|\xi|^2 + a_2(x) (1-e^{-i\xi}) + a_3(x) |\xi|$$
where
$$a_1(x)=\begin{cases}1 &, \text{ if }x<-6\\ 1-\frac{x+6}{2} &, \text{ if } x \in (-6,-4)\\ 0 &, \text{ otherwise } \end{cases},\ \ \ \ a_2(x)=\begin{cases}1 &, \text{ if }x\in[-4,4]\\ \frac{x+6}{2} &, \text{ if } x \in (-6,-4)\\ 1-\frac{x-4}{2} &, \text{ if } x\in(4,6)\\ 0 &, \text{ otherwise } \end{cases}$$
and 
$$a_3(x)=\begin{cases}1 &, \text{ if }x>6\\ \frac{x-4}{2} &, \text{ if } x \in (4,6)\\ 0 &, \text{ otherwise } \end{cases}.$$

\begin{figure}[!ht]
\begin{center}
\includegraphics[width=4in]{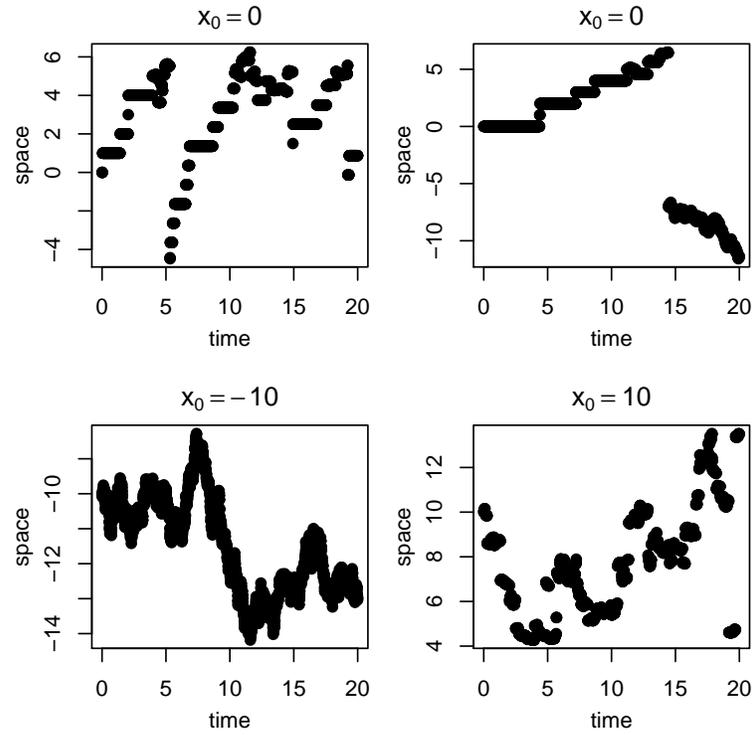}
\end{center}
\caption{
{\bf Brownian-Poisson-Cauchy-mixture Feller process.} Around the origin it behaves like a Poisson process, for smaller values like Brownian motion and a for larger values like a Cauchy process.
}
\label{fig:levy-mixture}
\end{figure}

These functions are Lipschitz continuous and thus a corresponding Feller process exists. Figure \ref{fig:levy-mixture} shows some samples of this process on $[0,20]$ with time-step size $\frac{1}{100}$. One can observe that the process behaves like a Poisson process around the origin, like a Cauchy process above 6 and like Brownian motion below -6.

\subsubsection*{Symmetric stable-like process}
A L\'evy process $L_t$ is a symmetric-$\alpha$-stable process if there exists an $\alpha\in(0,2]$ such that its characteristic function is given by
$$\E e^{i\xi L_t} = e^{-t |\xi|^\alpha}.$$
If we now define a function $x \mapsto \alpha(x)$ where $\alpha(x)$ takes only values in $(0,2]$ then there exists a family of of L\'evy processes $(L_t^{(x)})_{t\in[0,\infty)}$ such that for fixed $x$ the the process $L_t^{(x)}$ has the characteristic function
$$\E e^{iL_t^{(x)}\xi}  = e^{-t |\xi|^{\alpha(x)}}.$$
A corresponding Feller process exists and is unique if the function $x\mapsto \alpha(x)$ is Lipschitz continuous and bounded away from 0 and 2 \cite{Bass1988a}.

\begin{figure}[!ht]
\begin{center}
\includegraphics[width=4in]{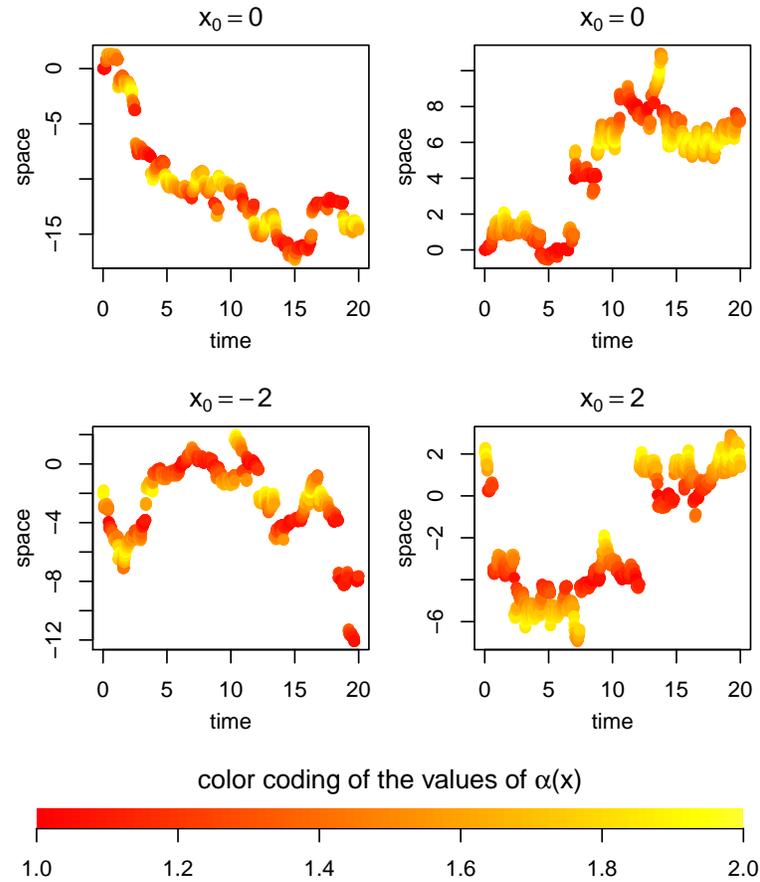}
\end{center}
	\caption{{\bf Stable-like processes with }$\mathbf{\alpha(x) := 1+\frac{19}{10}\left(\left(\frac{x}{4}-\left\lfloor \frac{x}{4}\right\rfloor\right) \land \left(\left\lceil \frac{x}{4}\right\rceil- \frac{x}{4}\right)\right)}${\bf.} Each position is color coded by the corresponding value of the exponent. In the yellow regions it behaves almost like Brownian motion, in the red regions it behaves almost like a Cauchy process.}
	\label{fig:stable-like}
\end{figure}

Figure \ref{fig:stable-like} shows some samples of a stable-like Feller process on $[0,20]$ with time-step size $\frac{1}{10}$ and  
$$\alpha(x) := 1+\frac{19}{10}\left(\left(\frac{x}{4}-\left\lfloor \frac{x}{4}\right\rfloor\right) \land \left(\left\lceil \frac{x}{4}\right\rceil- \frac{x}{4}\right)\right),$$
i.e.\ $x\mapsto a(x)$ is a function which is Lipschitz continuous (but not smooth) oscillating between 1 and (nearly) 2. To understand the figure note that we color coded the state space: red indicates $\alpha \approx 1$, yellow indicates $\alpha \approx 2$ and the values between these extremes are colored with the corresponding shade of orange. Now one can observe that the process behaves in the red areas like a Cauchy process and the more \textit{yellow} the state becomes, the more the process behaves like Brownian motion.

\subsubsection*{Normal inverse Gaussian-like process}
\label{ghp}
\begin{figure}[!ht]
\begin{center}
\includegraphics[width=4in]{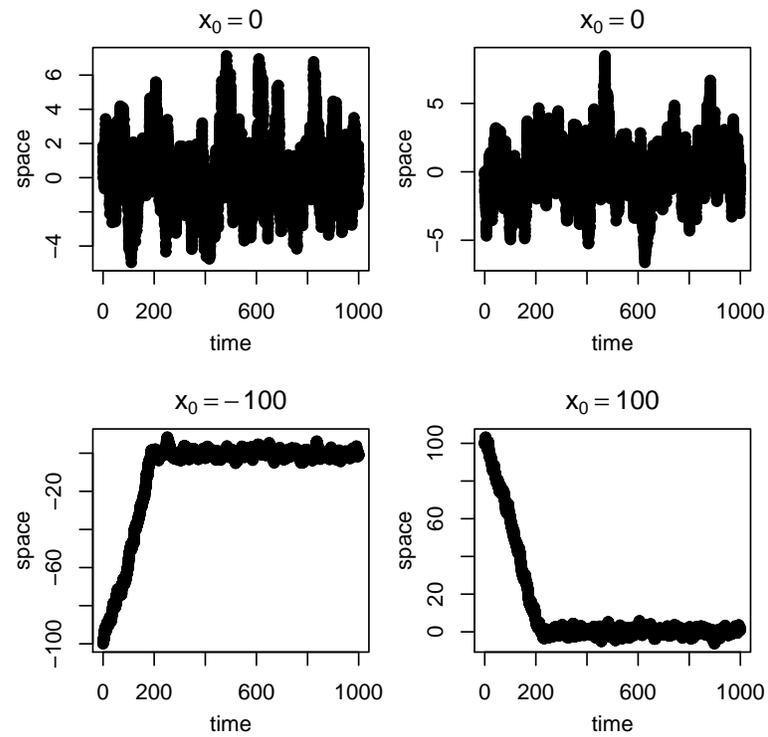}
\end{center}
	\caption{{\bf Normal inverse Gaussian-like processes with }$\mathbf{\beta(x):=-\frac{1}{\pi} \arctan(x)}${\bf.} The process features mean reversion to 0.}
	\label{fig:normal-inverse-gaussian-like}
\end{figure}

The characteristic function of a normal inverse Gaussian process $L_t$ is given by
$$\E e^{i\xi L_t} = \exp\left({t i\mu - t\delta \left(\sqrt{\alpha^2 - (\beta + i\xi)^2}- \sqrt{\alpha^2-\beta^2}\right)}\right)$$

where $\alpha> 0,\, -\alpha<\beta<\alpha,\, 0 < \delta$ and $\mu$. If we replace $\alpha,\, \beta,\, \delta,\, \mu$ by arbitrarily often differentiable bounded functions $\alpha(x),\, \beta(x),\, \delta(x),\, \mu(x)$ with bounded derivatives and assume that ther exist constants $c,C>0$ such that $\delta(x)>c,\, \alpha(x)-|\beta(x)| > c,\, c\leq \mu(x) \leq C$, then it was stated in \cite{BarnLeve2001} that a corresponding Feller process exists. Therein was also proposed an example of a mean reverting normal inverse Gaussian-like process, a special case of this model with mean 0 is obtained by setting
$$\alpha(x):=\delta(x):= 1,\ \mu(x):= 0\ \text{ and } \beta(x):=-\frac{1}{\pi} \arctan(x).$$
Note that the mean reversion is not introduced by using simply a drift which drags the process back to the origin. It is the choice of $\beta$ which yields an asymmetric distribution that moves the process back to the origin. The mean reversion can be observed in Figure \ref{fig:normal-inverse-gaussian-like} which shows samples of the normal inverse Gaussian-like process on $[0,1000]$ with time-step size $\frac{1}{10}$.

\subsubsection*{Meixner-like process}
\label{meix}
The characteristic function of a Meixner process $L_t$ is given by
$$\E e^{i\xi L_t} = \left(\frac{\cos \dfrac{b}{2}}{\cosh{\dfrac{a\xi-ib}{2}}}\right)^{2rt} \exp\left({imt\xi}\right)$$
where $a,r > 0,\, -\pi < b < \pi,\, m\in\R$. Details can be found in \cite{Grig1999}.

A family of Meixner processes which corresponds to a Feller process can be constructed by substituting the parameters $a,\, b,\, r,\, m$ by arbitrary often differentiable bounded functions $a(x),\, b(x),\, r(x),\, m(x)$ with bounded derivatives. The functions have to be bounded away from the critical values, i.e.\ $0 < a_0 \leq a(x),\, 0 < r_0 \leq r(x)$ and $-\pi < b_- < b(x) < b_+ < \pi$ for some fixed $r_0,\, a_0,\, b_-,\, b_+$. For further details see \cite{BoetJaco2003}. 

\begin{figure}[!ht]
\begin{center}
\includegraphics[width=4in]{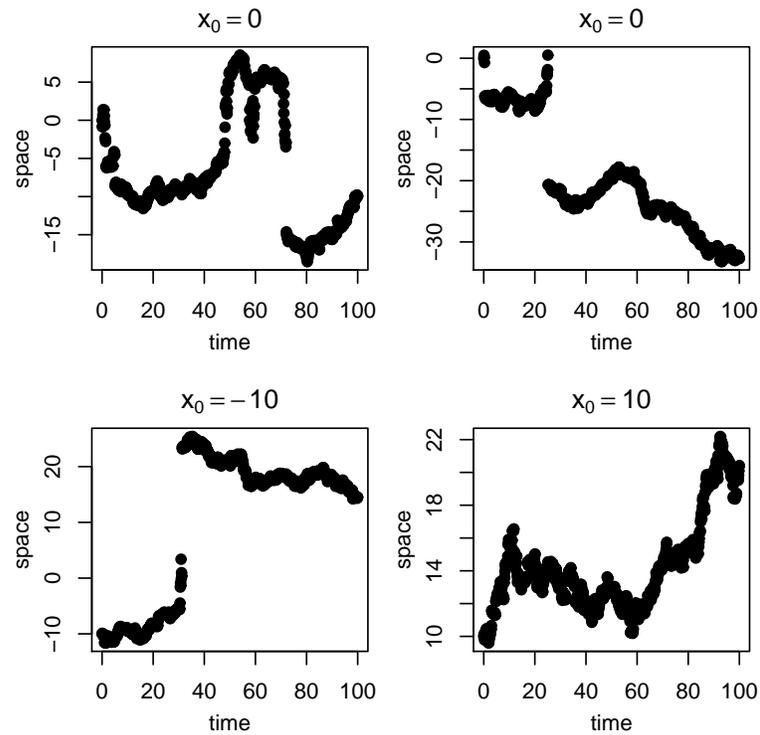}
\end{center}
\caption{
{\bf Meixner-like process with }$\mathbf{a(x):= 1+ 10 e^{-\frac{1}{25-x^2}}1_{(-5,5)}(x)}${\bf.}  The process moves \textit{with bigger steps} around the origin than for larger (and smaller) values. In fact by the choice of $x\mapsto a(x)$ the rate of the exponential decay of the transition density is reduced around the origin.
}
\label{fig:meixner-like}
\end{figure}

Figure \ref{fig:meixner-like} shows some samples of the Meixner-like process on $[0,100]$ with time-step size $\frac{1}{10}$ and 
$$b(x):=m(x):= 0,\ r(x):=1\ \text{ and } a(x):= \begin{cases}1+ 10 e^{-\frac{1}{25-x^2}} &, \text{ if } |x|< 5\\ 1 &, \text{ otherwise}\end{cases}.$$ 
The chosen functions satisfy the existence conditions from above. Furthermore the function $x\mapsto a(x)$ yields that the process \textit{moves with bigger steps} around the origin, to be precise: the Meixner distribution has semiheavy tails \cite{Grig2001a} and the parameter $a$ determines the rate of the exponential decay factor for the density. The effect on the sample path can be observed in Figure \ref{fig:meixner-like}.

\subsection*{Conclusion}
Using the presented mixture approach one can easily construct Feller models based on given L\'evy models. In these cases the existence of the process is granted. 

Furthermore the presented approximation is a very intuitive way to generate the sample path of a Feller processes. Obviously the method requires that one can simulate the increments of the corresponding L\'evy processes. But for L\'evy processes used in applications, especially together with Monte Carlo techniques, this poses no new restriction. 

Thus all necessary tools are available to use Feller processes as models for a wide range of applications.

\section*{Materials and Methods}

The simulations where done in R \cite{R2010} and the source code of the figures can be found as supporting information (Appendix S1).

\section*{Acknowledgments}


\begin{thebibliography}{10}

\bibitem{Baeu2005}
B.~Baeumer, M.~Meerschaert, and J.~Mortensen.
\newblock Space-time fractional derivative operators.
\newblock {\em Proc Amer Math Soc}, 133(8):2273--2282, 2005.

\bibitem{BarnLeve2001}
O.~Barndorff-Nielsen and S.~Levendorski{\u{\i}}.
\newblock Feller processes of normal inverse {G}aussian type.
\newblock {\em Quant Finance}, 1:318--331, 2001.

\bibitem{Bass1988a}
R.~Bass.
\newblock Uniqueness in law for pure jump {M}arkov processes.
\newblock {\em Probab. Theory and Relat. Fields}, 79:271--287, 1988.

\bibitem{Boet2005}
B.~B{\"o}ttcher.
\newblock A parametrix construction for the fundamental solution of the
  evolution equation associated with a pseudo-differential operator generating
  a {M}arkov process.
\newblock {\em Math. Nachr.}, 278(11):1235--1241, 2005.

\bibitem{Boet2008}
B.~B{\"o}ttcher.
\newblock Construction of time inhomogeneous {M}arkov processes via evolution
  equations using pseudo-differential operators.
\newblock {\em J. Lon. Math. Soc.}, 78(2):605--621, 2008.
\newblock to appear.

\bibitem{BoetJaco2003}
B.~B{\"o}ttcher and N.~Jacob.
\newblock Remarks on {M}eixner-type processes.
\newblock {\em Probabilistic Methods in Fluids (eds. I.M. Davies et al.)},
  pages 35--47, 2003.
\newblock World Scientific Press.

\bibitem{BoetSchi2009}
B.~B{\"o}ttcher and R.~Schilling.
\newblock {A}pproximation of {F}eller processes by {M}arkov chains with
  {L}\'evy increments.
\newblock {\em Stoch Dyn}, 9(1):71–--80, 2009.

\bibitem{BoetSchn2009}
B.~B\"ottcher and A.~Schnurr.
\newblock {The Euler scheme for Feller processes}.
\newblock 2010.

\bibitem{Cher1974}
P.~Chernoff.
\newblock {\em Product Semigroups, Nonlinear Semigroups and Addition of
  unbounded Operators}.
\newblock Number 140 in Mem. AMS, 1974.

\bibitem{ContTank2004}
R.~Cont and P.~Tankov.
\newblock {\em Financial modelling with jump processes}.
\newblock Chapman \& Hall/CRC, 2004.

\bibitem{CoxIngeRoss1985}
J.~C. Cox, J.~E. Ingersoll, and S.~A. Ross.
\newblock A theory of the term structure of interest rates.
\newblock {\em Econometrica}, 53(2):385--408, 1985.

\bibitem{Ditl1996}
P.~D. Ditlevsen, H.~Svensmark, and S.~Johnsen.
\newblock {Contrasting atmospheric and climate dynamics of the last-glacial and
  Holocene periods}.
\newblock {\em Nature}, 379:810--812, 1996.

\bibitem{DoneMaller2005a}
R.~A. Doney and R.~A. Maller.
\newblock {Cram\'er's Estimate for a Reflected L\'evy Process}.
\newblock {\em Ann Appl Probab}, 15(2):1445--1450, 2005.

\bibitem{Fell1951}
W.~Feller.
\newblock Diffusion processes in genetics.
\newblock In J.~Neyman, editor, {\em Proceedings of the Second Berkeley
  Symposium on Mathematical Statistics and Probability}, pages 227--246.
  University of California Press, 1951.

\bibitem{Gonc2003}
V.~Y. Gonchar, A.~V. Chechkin, E.~L. Sorokovoi, V.~V. Chechkin, L.~I.
  Grigor'eva, and E.~D. Volkov.
\newblock {Stable L\'evy distributions of the density and potential
  fluctuations in the edge plasma of the U-3M torsatron}.
\newblock {\em Plasma Physics Reports}, 29(5):380--390, 2003.

\bibitem{Grig1999}
B.~Grigelionis.
\newblock Processes of {M}eixner type.
\newblock {\em Lithuanian Math J}, 39(1):33--41, 1999.

\bibitem{Grig2001a}
B.~Grigelionis.
\newblock Generalized $z$-distributions and related stochastic processes.
\newblock {\em Lithuanian Math J}, 41(3):239--251, 2001.

\bibitem{Hoh98hab}
W.~Hoh.
\newblock {\em Pseudo differential operators generating {M}arkov processes}.
\newblock Habilitationsschrift. Universit{\"a}t Bielefeld, 1998.

\bibitem{IkedWata1989}
N.~Ikeda and S.~Watanabe.
\newblock {\em Stochastic Differential Equations and Diffusion Processes},
  volume~24 of {\em Math. Library}.
\newblock North-Holland, 2nd edition, 1989.

\bibitem{Jaco20012005}
N.~Jacob.
\newblock {\em Pseudo-differential operators and {M}arkov processes}, volume
  I-III.
\newblock Imperial College Press, 2001-2005.

\bibitem{JacoPotr2004}
N.~Jacob and A.~Potrykus.
\newblock {Roth's method applied to some pseudo-differential operators with
  bounded symbols. A case study.}
\newblock {\em Rend. Cir. Mat. Palermo (Ser. II)}, 76:45--57, 2005.

\bibitem{JacoShir2002}
J.~Jacod and A.~N. Shiryaev.
\newblock {\em Limit {T}heorems for {S}tochastic {P}rocesses}.
\newblock Springer, 2nd edition, 2002.

\bibitem{Jans1999}
H.~Janssen, K.~Oerding, F.~van Wijland, and H.~Hilhorst.
\newblock L\'evy-flight spreading of epidemic processes leading to percolating
  clusters.
\newblock {\em Eur. Phys. J. B}, 7:137--145, 1999.

\bibitem{Koch1989}
A.~N. Kochubei.
\newblock Parabolic pseudodifferential equations, hypersingular integrals and
  {M}arkov processes.
\newblock {\em Math. USSR Izvestija}, 33:233--259, 1989.

\bibitem{Kolo2000}
V.~N. Kolokoltsov.
\newblock Symmetric stable laws and stable-like jump-diffusions.
\newblock {\em Proc. London Math. Soc.}, 80:725--768, 2000.

\bibitem{MaRoeckZhan1998}
Z.-M. Ma, M.~R{\"o}ckner, and T.-S. Zhang.
\newblock Approximation of arbitrary {D}irichlet processes by {M}arkov chains.
\newblock {\em Ann. Inst. Henri Poincare}, 34(1):1--22, 1998.

\bibitem{Meer2002a}
M.~M. Meerschaert, D.~A. Benson, H.-P. Scheffler, and P.~Becker-Kern.
\newblock Governing equations and solutions of anomalous random walk limits.
\newblock {\em Physical Review E}, 66:060102, 2002.

\bibitem{ProtTala97}
P.~Protter and D.~Talay.
\newblock The {E}uler scheme for {L}{\'e}vy driven stochastic differential
  equations.
\newblock {\em Ann Probab}, 25(1):393--423, 1997.

\bibitem{R2010}
{R Development Core Team}.
\newblock R: A language and environment for statistical computing, 2010.
\newblock {ISBN} 3-900051-07-0.

\bibitem{Reev2008}
D.~Reeves, D.~Benson, M.~Meerschaert, and H.~Scheffler.
\newblock {Transport of Conservative Solutes in Simulated Fracture Networks 2.
  Ensemble Solute Transport and the Correspondence to Operator-Stable Limit
  Distributions}.
\newblock {\em Water Resources Research}, 44:W05410, 2008.

\bibitem{ReynFrye2007}
F.~M. Reynolds~AM.
\newblock Free-flight odor tracking in drosophila is consistent with an optimal
  intermittent scale-free search. 2(4): e354. doi:.
\newblock {\em PLoS ONE}, 2(4):e354, 2007.

\bibitem{Sato99}
K.~Sato.
\newblock {\em L\'evy Processes and Infinitely Divisible Distributions}.
\newblock Cambridge University Press, 1999.

\bibitem{Gore2004}
E.~Scalas, R.~Gorenflo, and F.~Mainardi.
\newblock Uncoupled continuous-time random walks: Solution and limiting
  behavior of the master equation.
\newblock {\em Physical Review E}, 69:011107, 2004.

\bibitem{SchiSchn09}
R.~Schilling and A.~Schnurr.
\newblock {The Symbol Associated with the Solution of a Stochastic Differential
  Equation}.
\newblock {\em Electron J Probab}, 15:1369--1393, 2010.
\newblock To appear in: El J Probab.

\bibitem{Stro1975}
D.~Stroock.
\newblock Diffusion processes associated with levy generators.
\newblock {\em Probab Theory Relat Fields}, 32:209--244, 1975.

\bibitem{Stro2003}
D.~W. Stroock.
\newblock {\em Markov Processes from K. It\^o's Perspective}.
\newblock Princeton University Press, 2003.

\bibitem{Woyc2001}
W.~A. Woyczy{\'n}ski.
\newblock L{\'e}vy processes in the physical sciences.
\newblock In {\em L{\'e}vy Process-Theory and Applications}, pages 241--266.
  Birkh{\"a}user, 2001.

\bibitem{Zhan2007}
Y.~Zhang, D.~Benson, M.~Meerschaert, and E.~M. LaBolle.
\newblock {Space-fractional advection-dispersion equations with variable
  parameters: Diverse formulas, numerical solutions, and application to the
  MADE-site data}.
\newblock {\em Water Resources Research}, 43:W05439, 2007.

\end{thebibliography}

\newpage
\section*{Figures}

\clearpage
\section*{Appendix S1}
\begin{verbatim}
################################################################################
# Simulations of Feller processes - B. Boettcher
################################################################################
# Code by B. Boettcher (bjoern.boettcher {at} tu-dresden.de) for the paper:
# Feller processes: The Next Generation in Modeling. Brownian Motion, Levy
# processes and beyond.
#
# Code executed within R (http://www.r-project.org/), version 2.11.1
#
# Note: We opted for clarity, rather than for the most efficient code. The 
# functions are written for an educated user, i.e. the functions do in general
# not check if the entered parameters have appropriate values.
#
# If you want to use this code to simulate Feller processes - instead of just
# reproducing the figures of the paper - please note:
#
# rfellerprocess(maxtime, timestepsize, startpoint, r1dstep) generates the 
# sample path and r1dstep is one of
# rlevymixture      - for a process given by the mixture approach 
# rstablefamily     - for a stable-like process
# rnigfamily        - for a normal inverse Gaussian-like process
# rmeixnerfamily    - for a Meixner-like process
# 
# The components of the mixture can be modified at the beginning of the section 
# "the Levy mixture approach".
# The families can be modified in the functions rstablefamily, rnigfamily and 
# rmeixnerfamily, respectively.
# 
# Examples:
# plot(rfellerprocess(10,0.001,0,rlevymixture))
# plot(rfellerprocess(10,0.001,0,rstablefamily))
# plot(rfellerprocess(1000,0.1,0,rnigfamily))
# plot(rfellerprocess(1000,0.1,0,rmeixnerfamily))
#
################################################################################

rfellerprocess <- function(maxtime,timestepsize,startpoint,r1dstep){
# generates the approximation to the path of the Feller process
#
# maxtime = time up to which the path will be simulated
# timestepsize = time-step size of the increments
# startpoint = starting point of the process
# r1dstep = a one dimensional function for one step generation 
#           (the function should accept two parameters: 
#            the time-step size and the current position)
#
# returns a list containing the time and position of the process as xy.coords

  n = ceiling(maxtime/timestepsize)
  v = numeric(n+1)
  v[1] = startpoint
  for ( i in 1:n) {
    v[i+1] = r1dstep(timestepsize,v[i]) + v[i]
  }
  invisible(xy.coords((0:n)* timestepsize,v,xlab="time",ylab="space"))
}


################################################################################
# the Levy mixture approach
################################################################################
# (here we work with global variables for convenience)
#

n = 3                                           # Number of processes to mix
rlevy = list(n)                                 # the processes
a = list(n)                                     # the coefficients
rlevy[[1]] = function(t) rnorm(1,sd=sqrt(t))    # L^1 = Brownian motion
rlevy[[2]] = function(t) rpois(1,lambda=t)      # L^2 = Poisson process
rlevy[[3]] = function(t) rcauchy(1,scale=t)     # L^3 = Cauchy process

# Now one could define just the functions directly (and the process exists
# always if the functions are non-negative, bounded and Lipschitz continuous)
# a[[1]] = 
# a[[2]] = 
#
# Instead, we define regions and generate the functions by linear interpolation

regions = c( 
-6,         # right endpoint of the region where the process is just L^1
-4,4,       # left & right endpoint of the region where the process is just L^2
6           # left endpoint of the region where the process is just L^3
)

# The following yields functions 'a' such that they are =1 on their region and 
# =0 on the others. Constructed by linear interpolation.
 
a[[1]] = approxfun(c(regions[1],regions[2]),c(1,0), method="linear",1,0)
for(i in 2:(n-1)) {
   a[[i]]=approxfun(c(regions[i*2-3],regions[i*2-2],regions[i*2-1],regions[i*2])
                    ,c(0,1,1,0),method="linear",0,0)
}
a[[n]] = approxfun(c(regions[n*2-3],regions[n*2-2]),c(0,1), method="linear",0,1)

rlevymixture = function(t,x) {
# generates an increment of the Levy family constructed as mixture with the 
# above defined processes and coefficients
#
# t = time-step size
# x = current position

  y = 0
  for(i in 1:n)
   y = y + a[[i]](x)*rlevy[[i]](t)
  return(y)
}


################################################################################
# the symmetric-alpha-stable family 
################################################################################

rstable = function(n,t,alpha=1) {
# generation of increments of an alpha stable process
#
# n = number of increments
# t = time-step size
# alpha = index of stability
#
# We use 
# * the algorithm form Chambers, J. M.; Mallows, C. L. & Stuck, B. W. 
#   A Method for Simulating Stable Random Variables, Journal of the American 
#   Statistical Association, 1976, 71, 340-344
# * the stability, i.e.: 
#   X has char.exp. |xi|^alpha ==> t^(1/alpha) X has char.exp. t |xi|^alpha

   W = rexp(n)
   Theta = runif(n,-pi/2,pi/2)
   return(t^(1/alpha)*(sin(alpha*Theta)/cos(Theta)^(1/alpha) 
           * (cos((1-alpha)*Theta)/W)^((1-alpha)/alpha)))
}

rstablefamily = function (t,x){
# generates an increment of the symmetric alpha stable process family
# 
# t = time-step size of the increment
# x = current position

  exponent = function(x) 1+19/10*pmin((x/4-floor(x/4)),(ceiling(x/4)-x/4))
                                                       # 0 < exponent < 2

  rstable(1,t,exponent(x))
}


################################################################################
# the normal inverse Gaussian family 
################################################################################

rinversegaussian = function(t,a,b) {
# generates an increment of the inverse Gaussian process 
#
# t = time-step size
# a,b = parameters of the inverse Gaussian process
#
# Algorithm from: V. Seshadri, The inverse Gaussian distribution: a case study 
# in exponential families, Oxford University Press, 1993 (page 203)

  y = rnorm(1)^2
  u = runif(1)
  m = t*a
  
  xtemp = m+m^2*y/(2*b)- m*sqrt(4*m*b*y+m^2*y^2)/(2*b)
  if (u <= m/(m+xtemp)) return(xtemp)
  else return(m^2/xtemp) 
}

rnormalinversegaussian = function(t,alpha,beta,delta) {
# generates an increment of the normal inverse Gaussian process
#
# t = time-step size
# alpha,beta,delta = parameters of the normal inverse Gaussian process
#
# Simulation based on subordination of Brownian motion by an inverse Gaussian
# process

  y = rinversegaussian(t,delta,sqrt(alpha^2-beta^2))
  return(beta*y+rnorm(1,sd=sqrt(y)))
}

rnigfamily = function(t,x) {
# generates an increment of the normal inverse Gaussian family
#
# t = time-step size
# x = current position

  alpha = function(x) 1                                #  alpha > 0
  beta = function(x) -1/pi*atan(x)                     # -alpha < beta < alpha
  delta = function(x) 1                                #  delta > 0
  mu = function(x) 0                                   # mu is a real number

  return(mu(x)*t+rnormalinversegaussian(t,alpha(x),beta(x),delta(x)))
}


################################################################################
# the Meixner family
################################################################################
# we use the package cxxPack for an implementation of the complex valued 
# Gamma function "cgamma()"
# to install the package use:
# install.packages("cxxPack")
#
# Please note, that we have no relation to the author of that package.
library("cxxPack") 

rmeixner = function(t,a,b=0,m,r) {
# generates an increment of the Meixner process
#
# t = time-step size
# a,b,m,r = parameters of the Meixner process
#
# The code is based on the method proposed in the preprint:
# Reiichiro Kawai, Parameter Sensitivity Estimation for Meixner Distribution and
# Levy Processes, University of Leicester
#
# We only implemented the method for b=0 and checked the results of the proposed
# method by an optical comparison of the resulting histogram with the known 
# density function for various parameters:
#
# dmeixner = function(t,a,b,m,r,x) {
#   density function of the increments of the Meixner process
#
#   r=r*t
#   m=m*t
#   return((2*cos(b/2))^(2*r)/(2*a*pi*gamma(2*r))
#           *exp(b*(x-m)/a)*abs(cgamma(r+1i*(x-m)/a))^2)
# } 
# 
# t =
# a =
# m =
# r =
#
# x=replicate(10000,rmeixner(t,a,0,m,r))
# hist(x,freq=F)
# f = function(x) dmeixner(t,a,0,m,r,x)
# curve(f,col="red",add=T)
# 

  if (b!=0) {
   print("Error: Only the case b=0 is implemented")
   return("Error")
  }
  r=t*r
  m=t*m
  repeat {
    Q = runif(1,min=-1,max=1)/runif(1,min=-1,max=1)
    V = a/2*max(sqrt(2*r),2*r)*Q
    U = runif(1)
    if (abs(Q) < 1) {
      if (gamma(r)^2*U < abs(cgamma(r+1i*V/a))^2) { 
        break
      }  
    } 
    else {
      if (!is.na(cgamma(r+1i*V/a))) 
     # Note: if abs(V) is to large then cgamma yields NaN (but in fact the value 
     # is very small, so also the following inequality would be false. So we 
     # just treated the NaN exception separately.) 
        if (max(1,2*r)*a^2*gamma(r+1)^2*U/(2*r) < abs(cgamma(r+1i*V/a))^2 *V^2){ 
             break
           }
    }     
  }
  return(V+m)
}

rmeixnerfamily = function(t,x) {
# generates an increment of the Meixner family
#
# t = time-step size
# x = current position

  a = function(x) if (abs(x)>=5) return(1) else return(1+10*exp(-1/(25-x^2)))
  # a>0
  b = function(x) 0                 # -pi < b < pi (implemented only for b=0 !!)
  m = function(x) 0                 # m is a real number
  r = function(x) 1                 # b>0
 
  return(rmeixner(t,a(x),b(x),m(x),r(x)))
}


################################################################################
################################################################################
# Generation of the Figures
################################################################################

# initialization of the random number generator and some graphic parameters
set.seed(1034)

par(oma = c(0,0,0,0), mar=c(4,3.2,2,0.8))
layout(matrix(c(1,2,3,4),nrow=2,byrow=TRUE),respect=TRUE)

plottype = "s"
plotpch = 19

figplot = function(xy,...) {
  plot(xy,xlab="",ylab="",pty=plottype,pch=plotpch,...)
  title(xlab="time",ylab="space",line=2.2)
}

# the levy mixture
figplot(rfellerprocess(20,0.01,0,rlevymixture),main=expression(x[0] == 0))
figplot(rfellerprocess(20,0.01,0,rlevymixture),main=expression(x[0] == 0))
figplot(rfellerprocess(20,0.01,-10,rlevymixture),main=expression(x[0] == -10))
figplot(rfellerprocess(20,0.01,10,rlevymixture),main=expression(x[0] == 10))
dev.copy2eps(file="Brownian-Poisson-Cauchy-mixture.eps",width = 4.8, height = 4.8)

# the stable-like process
layout(matrix(c(1,2,3,4,5,5),nrow=3,byrow=TRUE),heights=c(1,1,0.4), respect=TRUE)
par(oma = c(0,0,0,0), mar=c(4,3.2,2,0.8),cex=0.83)
colorexponent = function(x) pmin((x/4-floor(x/4)),(ceiling(x/4)-x/4))
cols = heat.colors(100)

path=rfellerprocess(20,0.01,0,rstablefamily)
figplot(path,main=expression(x[0] == 0)
  ,col=cols[1+ceiling(colorexponent(path$y)*150)])
path=rfellerprocess(20,0.01,0,rstablefamily)
figplot(path,main=expression(x[0] == 0)
  ,col=cols[1+ceiling(colorexponent(path$y)*150)])
path=rfellerprocess(20,0.01,-2,rstablefamily)
figplot(path,main=expression(x[0] == -2)
  ,col=cols[1+ceiling(colorexponent(path$y)*150)])
path=rfellerprocess(20,0.01,2,rstablefamily)
figplot(path,main=expression(x[0] == 2)
  ,col=cols[1+ceiling(colorexponent(path$y)*150)])

# legend for the index of the stable-like process
plotlegend = function() {
  par(mar=c(3,0,2,0))
  plot(c(1,2),c(0,1),type="n",axes=0,xlab="",ylab="")
  title(expression(paste("color coding of the values of ", alpha, "(",x,")")))
  axis(1)
  n=80
  xr = seq(1,n)/n+1
  xl = xr-1/n
  yb = rep(0,n)
  yt = rep(1,n)
  rect(xl,yb,xr,yt,col=cols[1:n],border= NA) 
}
plotlegend()
dev.copy2eps(file="stable-like.eps",width = 4.8, height = 5.76)

# the normal inverse Gaussian-like process
par(oma = c(0,0,0,0), mar=c(4,3.2,2,0.8))
layout(matrix(c(1,2,3,4),nrow=2,byrow=TRUE),respect=TRUE)

figplot(rfellerprocess(1000,0.1,0,rnigfamily),main=expression(x[0] == 0))
figplot(rfellerprocess(1000,0.1,0,rnigfamily),main=expression(x[0] == 0))
figplot(rfellerprocess(1000,0.1,-100,rnigfamily),main=expression(x[0] == -100))
figplot(rfellerprocess(1000,0.1,100,rnigfamily),main=expression(x[0] == 100))
dev.copy2eps(file="normal-inverse-gaussian-like.eps",width = 4.8, height = 4.8)

# the Meixner-like process
figplot(rfellerprocess(100,0.1,0,rmeixnerfamily),main=expression(x[0] == 0))
figplot(rfellerprocess(100,0.1,0,rmeixnerfamily),main=expression(x[0] == 0))
figplot(rfellerprocess(100,0.1,-10,rmeixnerfamily),main=expression(x[0] == -10))
figplot(rfellerprocess(100,0.1,10,rmeixnerfamily),main=expression(x[0] == 10))
dev.copy2eps(file="meixner-like.eps",width = 4.8, height = 4.8)

\end{verbatim}

\end{document}